\documentclass[12pt,reqno]{article}
\usepackage{amsmath,amsthm,amsfonts,amssymb,amscd}
\usepackage{graphicx}
\usepackage{amssymb}
\usepackage{amsmath}
\usepackage{amscd}
\newtheorem{theorem}{Theorem}

\newtheorem{definition}[theorem]{Definition}

\newtheorem{lemma}[theorem]{Lemma}

\newtheorem{remark}[theorem]{Remark}

\newcommand{\curl}{\operatorname{{curl}}}
\newcommand{\nat}{\mathbb{N}}

\newcommand\grad{\operatorname{\nabla}}
\newcommand\lame{\operatorname{\mathcal{L}}}
\newcommand\lamet{\operatorname{\tilde{\mathcal{L}}}}
\newcommand\V[1]{\mathbf{#1}}
\renewcommand{\div}{\operatorname{div}}
\newcommand{\parag}{\hspace*{0.5cm}}

\title{Asymptotic stability and regularity of solutions for a magnetoelastic system in
  bounded domains}
\author{J\'auber C. Oliveira\\
              Departamento de Matem\'atica,\\
              Universidade Federal de Santa Catarina,\\
              Florian\'opolis, SC, 88040-900, Brasil,\\
              e-mail: j.c.oliveira@ufsc.br 
}
\begin{document}
\maketitle
\begin{abstract}
We prove the existence of strong time-periodic solutions and their asymptotic stability with the total energy of the perturbations decaying to zero at an exponential decay rate as $t \rightarrow \infty$ for a semilinear (nonlinearly coupled) magnetoelastic system in bounded, simply connected,
three-dimensional domain. The mathematical model includes a mechanical dissipation and a periodic forcing
function of period $T$. In the second part of the paper, we consider a magnetoelastic system in the form of a semilinear initial
boundary value problem in a bounded, simply-connected two-dimensional domain. We use LaSalle invariance principle to obtain results on the asymptotic behavior of solutions. This second result was obtained for the system under the action of only one dissipation (the natural dissipation of the system).
keywords: magnetoelastic system, asymptotic stability, periodic solutions
MSC2010: 35B10, 35B40, 35B35, 35L71, 74F15\footnote{The author gratefully acknowledges the financial support of CNPq \\(Proc. 201247/2010-0).}
\end{abstract}

\section{Part I: Regularity and Stability for Solutions of a Magnetoelastic System in 3D}
\subsection{Introduction}
We consider a system of partial differential equations that models the interaction between an elastic body and a constant magnetic field $\tilde{\mathbf{H}}$ acting on it. The elastic body is
assumed to be conducting, non-ferromagnetic, homogeneous and isotropic. The body
occupies a bounded, simply-connected domain $\Omega \subset \mathbb{R}^3$ with
boundary of class $C^2$. Let $T$ be a fixed positive number and let
$\V{u}(\V{x},t)$ denote the displacement field. Let $\mathbf{f}(\mathbf{x},t)$ denote an external force
acting on the body. We describe the coupling between the elastic and the magnetic systems as follows. The displacement in the body induces a magnetic field within the body given by $\tilde{\V{H}}+\V{h}(\V{x},t)$. The magnetic field exerts a Lorentz force on the elastic body to change the displacement field.The corresponding system of partial differential equations (\cite{ErMau90} where $\V{\rho}=0$) is:
\begin{align}
&\rho_{_M} \V{u}^{\prime\prime} - \mu\operatorname{\Delta} \V{u} - (\lambda+\mu)\grad \div\V{u}
  +\V{\rho}(\V{u}^ {\prime})= \mu_0 (\curl\V{h}) \times \left(\V{h}+\tilde{\V{H}}\right)
  + \V{f},\label{eq_1.1}\\ 
&\V{h}^{\prime} + \nu_1\;\curl\curl \V{h} = \curl \left[\V{u}^{\prime} \times
    \left(\V{h}+\tilde{\V{H}}\right)\right], \label{eq_1.2}\\
&\div \V{h} = 0,\label{eq_1.3}
\end{align}
in $Q_T=\Omega \times (0,T)$  ($Q_T=\Omega \times \frac{\mathbb{R}}{T
  \mathbb{Z}}$ when we consider time-periodic solutions). $\V{u}^{\prime}$ is
the partial derivative of $\V{u}(\V{x},t)$ with respect to $t$ (time) and
$\V{u}^{\prime\prime}$ means $(\V{u}^{\prime})^{\prime}$. Vector fields are represented in bold
face. $\rho_{_M}$ is the mass density per unity volume of
the medium. $\lambda$ and $\mu$ are Lam\'e's constants $(\mu \geq 0$, $3
\lambda + 2 \mu \geq 0)$. $\nu_1 = 1/(\sigma \mu_0)$ where $\sigma>0$
represents the conductivity of the material and $\mu_0$ is a positive number
representing the magnetic permeability. $\V{\rho}(\V{u}^{\prime})$ represents a
dissipation that acts on the elastic body.

The initial conditions associated with this system are
\begin{equation}\label{eq_1.4A}
\V{u}(\V{x},0) = \V{u}_0(\V{x}), \quad \V{u}^{\prime}(\V{x},0) =
\V{u}_1(\V{x}),\quad \V{h}(\V{x},0) = \V{h}_0(\V{x}),
\end{equation}
while in the time-periodic setting, the periodicity conditions associated with
this system are
\begin{equation}\label{eq_1.4}
\V{u}(\V{x},0) = \V{u}(\V{x},T), \quad \V{u}^{\prime}(\V{x},0) =
\V{u}^{\prime}(\V{x},T),\quad \V{h}(\V{x},0) = \V{h}(\V{x},T).
\end{equation}
The boundary conditions are 
\begin{equation}\label{eq_1.5}
\V{u}=0, \quad \V{h} \cdot \V{n} = 0, \quad (\operatorname{curl} \V{h}) \times
\V{n} = 0\;\;\mbox{on}\;\;\Sigma_T=\partial \Omega \times (0,T),
\end{equation}
( in the time-periodic case: $\Sigma_T=\partial \Omega \times
\frac{\mathbb{R}}{T \mathbb{Z}}$) where $\V{n} = \V{n}(\V{x})$ denotes the
outer unit normal at $\V{x}=(x_1,x_2,x_3) \in \partial\Omega$.

System \eqref{eq_1.1}--\eqref{eq_1.3} consists of a second order hyperbolic
problem coupled with a parabolic equation similar to the Navier-Stokes
equations through a semilinear coupling.

The first investigation on the existence of solutions (with $\V{\rho}=0$) for
the magnetoelastic system in bounded domains was carried out by Botsenyuk
\cite{Bot92}, who proved the existence of weak solutions for the
initial-boundary value problem in $\Omega \times (0, T)$ under the initial and
boundary conditions \eqref{eq_1.4A} and \eqref{eq_1.5}. In \cite{Bot96a},
Botsenyuk proved the existence and uniqueness of a global strong solution under
the hypothesis that the data is sufficiently small.

Andreou and Dassios (\cite{AnDa97}) studied the Cauchy problem in $\mathbb{R}^3$
and showed, under some regular conditions on the initial data, that the solution
of the system decays to zero, at a polynomial rate, as time goes to infinity.

Menzala and Zuazua (\cite{MeZu98}) proved the asymptotic stability for the system
with linearized coupling and no mechanical dissipation (i.e., the only
dissipation acting was the natural dissipation of the system). The proof used
LaSalle's invariance principle and did not provide a decay rate for the
system. 

Rivera and Santos (\cite{RivSan03}) proved for domains of special type that the energy associated to the linear three-dimensional magnetoelastic system decays polynomially to zero as time goes to infinity, provided the initial data is smooth enough.

Char\~ao, Oliveira and Menzala (\cite{ChOlMe09}) proved that the total
energy of this system tends to zero as $t \rightarrow \infty$ when a nonlinear
dissipation $\V{\rho}(\V{x},\V{u}^{\prime})$ is effective on a small subregion
of the domain. The rate of decay is given and depends on the behavior of
$\V{\rho}$ with respect to the second variable: algebraic decay for strong
solutions and exponential decay if the behavior is close to linear for weak
solutions.

Mohebbi and Oliveira (\cite{MohOliv12}) proved the existence of weak time-periodic solutions for this magnetoelastic system under the effect of a nonlinear dissipation $\V{\rho}(\V{u}^{'})$ like $|\V{u}|_E^p \mathbf{u}$ with $p \in [3,4]$ (or the simplest case $p=0$). The solutions have same period as the given time-periodic  body force $\mathbf{f}(x,t)$. In that paper it is also argued that the existence of weak solutions holds for $p \geq 2$ (or $p=0$) in the two-dimensional case.

In this first part of the paper, we consider the regularity and stability of time-periodic solutions for the semilinear magnetoelastic system obtained by \cite{MohOliv12} under the
effect of a linear dissipation $\V{\rho}(\V{u}^{\prime})=\alpha \V{u}^{\prime}$ or a nonlinear dissipation
satisfying certain restrictive condition.This results were motivated by our desire to show asymptotic stability of time-periodic solutions for the magnetoelastic system.
\subsection{Function spaces}
\label{sec:FS}
\hspace{5mm} We consider a bounded, simply-connected domain $\Omega$ of
$\mathbb{R}^3$ with $\partial \Omega$ a 2-manifold of class
$C^2$. $C^\infty_0(\Omega)$ and $C^\infty_0(\overline\Omega)$ are the space of
smooth functions with compact support in $\Omega$ and the space of smooth
functions with compact support in $\mathbb{R}^3$ restricted to $\Omega$.
respectively. The usual Sobolev spaces are denoted by $W^{s,p}(\Omega),
\;W_{0}^{s,p}(\Omega)$, $s \geq 0$, with the simplification $H^s(\Omega) :=
W^{s,2}(\Omega),\;H_0^s(\Omega) := W_0^{s,2}(\Omega)$. The norm of $H^{s}$
is denoted by $\|\;\|_s$ and $H^{-s}(\Omega)$ indicates the dual space of
$H^s_0(\Omega)$. More details on some definitions that follow can be found in
\cite{DaLi90}.

We recall that
\[
H(\div,\Omega) := \{\V v \in L^2(\Omega):\; \div\V v \in L^2(\Omega)\},
\]
is a Hilbert space and $C_0^\infty(\overline\Omega)$ is dense in this space. The
trace map $\gamma_n$ defined on $C_0^\infty(\overline\Omega)$ by $\gamma_n(\V
v) = \V v \cdot \V n|_{\partial\Omega}$ can be extended to a linear continuous
mapping $\gamma_n$ from $H(\div,\Omega)$ onto $H^{-1/2}(\partial\Omega)$ which
maps $\V v$ into $\V v\cdot\V n|_{\partial\Omega}$. 
Let $$\mathcal{V}:=\left\{\V{h} \in C_{0}^{\infty}(\overline{\Omega}):
\div\V{h}=0,\;\V{h} \cdot \V n|_{\partial\Omega} = 0\right\},$$
$H^1_{\sigma}(\Omega)$ is the closure of $\mathcal{V}$ in $H^1(\Omega)$, which
can be characterized as $$H^1_{\sigma}(\Omega)=\left\{\V{h} \in
H^1(\Omega): \div\V{h}=0,\;\V{h} \cdot \V n|_{\partial\Omega} =
0\right\}.$$
Also, $L^2_{\sigma}(\Omega)$ is the closure of $\mathcal{V}$ in $L^2(\Omega)$,
which can be characterized as $$L^2_{\sigma}(\Omega)=\left\{\V{h} \in
L^2(\Omega): \div\V{h}=0,\;\V{h} \cdot \V n|_{\partial\Omega} =
0\right\}.$$
The Hilbert space $L^2_{\sigma}(\Omega)$ is equipped with the usual
$L^2$-norm, $|\;\;|_2$, and $L^2$-inner product, $(\;,\;)_2$; while
$H^1_{\sigma}(\Omega)$ is also a Hilbert space with inner product
\[
((\V{h},\V{h})):=(\curl \V{h},\curl\V{h})_{2}
\]
and norm $\|\;\|$ induced from this inner product (which is equivalent to the
usual $H^1$-norm). For $s \geq 0$, we define $H^{s}_{\sigma}(\Omega)$ as the
closure of $\mathcal{V}$ in $H^{s}(\Omega)$. $H^{s}_{\sigma}(\Omega)$ is a
Hilbert space equipped with the inner product of $H^{s}(\Omega)$, denoted by
$((\;,\;))_{s}$. As usual, $H^{-s}_\sigma(\Omega)$ indicates the dual of
$H^s_\sigma(\Omega)$.

Let $H(\curl,\Omega):=\{\V{v} \in L^2(\Omega):\;\curl \V{v} \in L^2(\Omega)\}$,
which is a Hilbert space with the norm
$\|\V{v}\|_{H(\curl,\Omega)}:= \{|\V{v}|_{2}^2 + |\curl \V{v}|_{2}^2\}^{1/2}$.
$C_0^\infty(\overline\Omega)$ is dense in $H(\curl;\Omega)$. The trace
$\gamma_{\tau}$ defined on $C_0^\infty(\overline\Omega)$ by
$\gamma_{\tau}(\V{v})=\V{v}\times \V n|_{\partial \Omega}$ can be extended
continuously to a linear continuous mapping $\gamma_{\tau}$ from
$H(\curl,\Omega)$ into $H^{-1/2}(\partial \Omega)$. We also
define $$U_2=\{\V{w} \in L^2_{\sigma} \cap H^2(\Omega): \curl \V{w} \in
H_0(\curl,\Omega)\},$$ equipped with norm defined by
$\|\V{w}\|_{U_2}=\left\{\|\V{w}\|_{H(\curl,\Omega)}^2 + |\curl \curl
\V{w}|_2^2\right\}^{1/2}$.

\subsection{Theorem on the existence of periodic solutions}
\label{sec:statement}
\hspace{5mm} For simplicity in notation, we consider the existence of periodic
solutions of period $T>0$ for the semilinear coupled system rewritten in the
following form:
\begin{align}\label{eq_0}
  & \V{u}^ {\prime\prime} + \lame \V{u}+\V{\rho}(\V{u}^ {\prime}) = \curl
  \V{h} \times \left(\V{h}+\tilde{\V{H}}\right) +
  \V{f},\nonumber\\
& \V{h}^{\prime} + \nu_1\; \lamet \V{h} = \curl\left[\V{u}^{\prime} \times
    \left(\V{h}+\tilde{\V{H}}\right)\right],
\end{align}
where
\begin{align*}
&\lame(\;) = -\mu \operatorname{\Delta}(\;) - (\lambda+\mu)\grad\div(\;),
  & &\mathcal{D}(\lame)=H^2(\Omega) \cap H_{0}^{1}(\Omega), \\
&\lamet(\;) = \curl\curl(\;),&  &\mathcal{D}(\lamet)=H^2(\Omega) \cap
  H_{\sigma}^{1}(\Omega).
\end{align*}

We will find a weak $T$-periodic solution for this system under the following
assumptions:\\
(H0) $\Omega \subset \mathbb{R}^3$ is a bounded, simply-connected domain with boundary $\partial \Omega$ of class $C^2$.\\
\noindent For $p \in [3,4]$, $q =(p+2)/(p+1)$:\\
(H1) $\V{f}$  belongs to $C([0,T],L^{q}(\Omega))$ and $\V{f}(0)=\V{f}(T)$.\\
(H2) The continuous function $\V{\rho}:\mathbb{R}^{3}\rightarrow \mathbb{R}^{3}$ satisfies:\\
There exists a positive constant $K_0$ such that $$(\V{\rho}(\V z),\V z)_{E} \geq
K_0 |\V z|_{E}^{p+2}, \qquad\forall \V z \in \mathbb{R}^{3},$$
There exist positive constants $r_{\rho}$ and $K_1$ such that
$$|\V{\rho}(\V z)|_{E} \leq K_{1} |\V z|_{E}^{p+1} \quad\text{if }|\V z|_{E}
  \geq r_{\rho}.$$
($|\;\;|_{E}$ and $(\;,\;)_{E}$ denote the Euclidean norm and corresponding
  inner product.)

Now, we state the theorem of existence of weak time-periodic solution proved by Mohebbi and Oliveira (2012) (\cite{MohOliv12}):
\begin{theorem}\label{thm_1}
Let $T>0$ be the period of the function $\V{f}$. We assume
\textup{(H0)--(H2)}. Then, there exist $\V{u} \in L^{\infty}(0,T;
H_{0}^{1}(\Omega)$ with $\V{u}^{\prime} \in L^{\infty}(0,T; L^{2}(\Omega)) \cap
L^{p+2}(0,T;L^{p+2}(\Omega))$ and $\V{h} \in
L^{\infty}(0,T;L^2_{\sigma}(\Omega)) \cap L^2(0,T;H_{\sigma}^{1}(\Omega))$ such
that $(\V{u},\V{h})$ satisfies the following system:
\begin{align}
& -\int_{0}^{T} \left(\V{u}^{\prime},\V{\varphi}\right)_{2} \eta^{\prime} ds
  +\int_{0}^{T}a_{I}\left(\V{u},\V{\varphi}\right) \eta ds +\int_{0}^{T}
  \left(\V{\rho}(\V{u}^{\prime}),\V{\varphi}\right)_{2} \eta
  ds \nonumber\\
&\qquad\qquad\qquad\qquad\quad=
\int_{0}^{T}B_{I}(\V{h}, \V{h})(\V{\varphi}) \eta\;ds
 +\int_{0}^{T} \left(\V{f},\V{\varphi}\right)_{2}\;{\eta} ds, \nonumber\\
&\qquad\qquad\qquad\qquad\qquad\qquad\qquad\quad\forall \V{\varphi} \in
 H_{0}^{1}(\Omega)\cap L^{p+2}(\Omega),\quad\forall \eta \in
 \mathbb{D}_{T},\\
& -\int_{0}^{T} \left(\V{h},\V{\psi}\right)_{2}\;\tilde{\eta}^{\prime}\;ds
  +\int_{0}^{T} a_{II}(\V{h},\V{\psi}) ds =
  \int_{0}^{T}B_{II}(\V{u}^{\prime},\V{h})(\V{\psi})\tilde{\eta} ds,\nonumber\\
& \qquad\qquad\qquad\qquad\qquad\qquad\qquad\qquad\qquad\quad
  \forall \V{\psi} \in H^{3/2}_{\sigma}(\Omega),\quad \forall \tilde{\eta} \in
  \mathbb{D}_{T},\label{eq_2.2}
\end{align}
where $$\mathbb{D}_{T}:=\left\{\omega \in
C^{\infty}(\mathbb{R}):\;\omega(s)=\omega(s+T),\;\forall s \in
\mathbb{R}\right\},$$
$$a_{II}(\V{h},\V{w}):=\nu_1\;((\V{h},\V{w} )), \;\forall \V{h},\V{w} \in
H^1_{\sigma}(\Omega),$$ 
$$a_{I}(\V{u},\V{v}):=\mu\;(\grad \V{u}, \grad \V{v})_2 + (\lambda+\mu) (\div
\V{u},\div\V{v})_2,\;\forall \V{u},\V{v} \in H_0^1(\Omega),$$
$$B_{II}(\V{w},\V{b})(\;):=\int_{\Omega}\V{w} \times
  \left(\V{b}+\tilde{\V{H}}\right) \cdot \curl(\;) d\V{x}$$
and
$$B_{I}(\V{h},\V{b})(\;):=\int_{\Omega} (\curl \V{h}) \times
\left(\V{b}+\tilde{\V{H}}\right) \cdot (\;) d\V{x},$$
with $B_{II}: L^2(\Omega) \times H_{\sigma}^{1}(\Omega) \rightarrow
H^{-3/2}_{\sigma}(\Omega)$ and $B_{I}:H_{\sigma}^{1}(\Omega) \times
H_{\sigma}^{1/2}(\Omega) \rightarrow H^{-1}(\Omega)$. Furthermore, $\V{u} \in
C(0,T;H^{1/2}(\Omega))$ and $\V{h} \in C(0,T;H_{\sigma}^{-1/4}(\Omega))$.
\end{theorem}
\subsection{Strong Solution for the Initial Boundary-Value Problem without Mechanical Dissipation}
We consider now the initial-boundary value problem (IBVP) described by the following semilinear magnetoelastic system:
\begin{align}\label{goveq}
& \mathbf{u}^{''} + \mathcal{L} \mathbf{u} =\curl \mathbf{h} \times \left(\mathbf{h}+\tilde{\V{H}}\right) + \mathbf{f},\nonumber\\
& \mathbf{h}^{'} + \nu_1 \tilde{\mathcal{L}} \mathbf{h} = \curl\left[\mathbf{u}^{'} \times \left(\mathbf{h}+\tilde{\V{H}}\right)\right],\\
& \mbox{div}\;\mathbf{h}=0.\nonumber
\end{align}
\subsection{A Priori Estimates}
We consider a priori estimates for solutions of the equations in the Faedo-Galerkin approach. For this we consider a basis $\{\V{b}_j\}_{j \in \mathbb{N}}$ of eigenfunctions of
 $$a_{I}(\V{b}_j, \V{\varphi})=\lambda_j (\V{b}_j,\V{\varphi})_2,\forall \V{\varphi} \in H^{1}_{0}(\Omega)$$ and a basis $\{\tilde{\V{b}}_k\}_{k \in \mathbb{N}}$ of eigenfunctions (\cite{GeBrLe06}, p. 58) of 
 $$a_{II}(\tilde{\V{b}}_k,\psi)=\tilde{\lambda}_k (\tilde{\V{b}}_k,\V{\psi})_2,\forall \V{\psi}  \in H^{1}_{\sigma}(\Omega).$$ 
 
Given $m \in \nat$, we define the finite-dimensional spaces $$S_m=\text{span}\;\{\V{b}_1,\V{b}_2,\cdots,\V{b}_m\}\quad \text{and}\quad \tilde{S}_m=\text{span}\;\{\tilde{\V{b}}_1,\tilde{\V{b}}_2,\cdots,\tilde{\V{b}}_m\},$$ so that $$\V{u}_m = \sum_{j=1}^{m} c_j b_j,\quad \V{h}_m=\sum_{j=1}^{m}\tilde{c}_j \tilde{\V{b}}_j.$$ We consider the weak formulation written in terms of $\{u_m\}$ and $\{h_m\}$ in the standard way.

Let us define the energy of the system as
\begin{align}
2 \mathcal{E}(t) := |\mathbf{u}_m^{'}(t)|_2^2 + \|\mathbf{u}_m(t)\|^2 + |\mathbf{h}_m(t)|_2^2.\nonumber
\end{align}

It follows from the weak formulation with $\V{\varphi}=\mathbf{u}_m^{'}$ and $\V{\psi}= \mathbf{h}_m$, adding the resulting equations, the basic energy estimate
\begin{align}\label{ineq0}
\frac{d\mathcal{E}}{dt}+  \nu_1 \|\mathbf{h}_m\|^2& \leq  |\mathbf{f}|_2\;|\mathbf{u}_m^{'}|_2
\end{align}
This estimate implies that $$\sqrt{\mathcal{E}(t)} \leq \sqrt{\mathcal{E}(0)}+ \|\V{f}\|_{L^{1}(0,T;L^2(\Omega))}, \forall t \in [0,T]$$ using arguments similar to the ones that we will show for the second energy.

It also follows from the weak formulation with $\V{\varphi}=\mathcal{L}\mathbf{u}_m^{'}$ and $\V{\psi}=\tilde{\mathcal{L}} \mathbf{h}_m$, adding the resulting equations, the following estimate
\begin{align}\label{ineq1}
\frac{d}{dt}\left\{\|\mathbf{u}_m^{'}\|^2  + |\mathcal{L}\mathbf{u}_m|^2 + \|\mathbf{h}_m\|^2\right\} + 2 \nu_1 |\tilde{\mathcal{L}}\mathbf{h}_m|^2& \leq 2 c_{\mu}\;|\tilde{\mathcal{L}}\mathbf{h}_m|^2 \|\mathbf{u}_m^{'}\|\nonumber\\
& + 2 \|\mathbf{f}\|\;\|\mathbf{u}_m^{'}\|
\end{align}
where $c_{\mu}$ is a positive constant.\\
Therefore,
\begin{align}\label{ineq2}
\mathcal{E}_1(t) +\nu_1 \int_{0}^{t} |\tilde{\mathcal{L}}\mathbf{h}_m|^2 ds & \leq \mathcal{E}_1(0) + \sqrt{2} c_{\mu} \int_{0}^{t} \mathcal{E}_1^{1/2}(s) |\tilde{\mathcal{L}}\mathbf{h}_m|^2 ds\nonumber\\
& + \sqrt{2} \;\int_{0}^{t}\|\mathbf{f}(s)\|\;\mathcal{E}_1(s)^{1/2}\; ds
\end{align}
where
\begin{align}
 2 \mathcal{E}_1(t) := \|\mathbf{u}_m^{'}(t)\|^2 + |\mathcal{L}\mathbf{u}_m(t)|^2 + \|\mathbf{h}_m(t)\|^2.
\end{align}
Brezis inequality implies that
\begin{align}\label{ineq3}
\mathcal{E}_1^{1/2}(t) &\leq \mathcal{E}_1^{1/2}(0) +
 \sqrt{2} c_{\mu} \int_{0}^{t} |\tilde{\mathcal{L}}\mathbf{h}_m|^2 ds + \sqrt{2} \;\int_{0}^{t}\|\mathbf{f}(s)\|\; ds\nonumber\\
& \leq c_0 + c_1 \int_{0}^{t} |\tilde{\mathcal{L}}\mathbf{h}_m|^2 ds + c_2 \int_{0}^{T}\|\mathbf{f}(s)\|\; ds
\end{align}
Returning to inequality \eqref{ineq2}, we get
\begin{align}
\nu_1 \int_{0}^{t} |\tilde{\mathcal{L}}\mathbf{h}_m|^2 ds & \leq c_0^2 +
c_0 c_1 \int_{0}^{t} |\tilde{\mathcal{L}}\mathbf{h}_m|^2 ds + c_1^2 \left(\int_{0}^{t} |\tilde{\mathcal{L}}\mathbf{h}_m|^2 ds \right)^2\nonumber\\
& + 2 c_1 c_2 \left(\int_{0}^{T}\|\mathbf{f}(s)\|\; ds\right) \left(\int_{0}^{t} |\tilde{\mathcal{L}}\mathbf{h}_m|^2 ds\right)\nonumber\\
&  + c_0 c_2 \int_{0}^{T}\|\mathbf{f}(s)\| ds  +  c_2^2 \left(\int_{0}^{T}\|\mathbf{f}(s)\|\; ds\right)^2.
\end{align}
Let $Z:=  \int_{0}^{t} |\tilde{\mathcal{L}}\mathbf{h}_m|^2 ds$. Then, we must have
\begin{align}
0 & \leq c_1^2 Z^2 + \left\{c_0 c_1 + 2 c_1 c_2 \left(\int_{0}^{T}\|\mathbf{f}(s)\|\; ds \right) -\nu_1\right\} Z \nonumber\\
& + c_0 c_2 \int_{0}^{t}\|\mathbf{f}(s)\| ds  +  c_2^2 \left(\int_{0}^{t}\|\mathbf{f}(s)\|\; ds\right)^2 + c_0^2.
\end{align}
We recall also an auxiliary lemma due to Botsenyuk's (\cite{Bot96a}).
\begin{lemma}
Assume that $x,\gamma \in C^{0}([0,T];\mathbb{R})$ are non-negative functions satisfying the inequality
\begin{align}
x(t) \leq \gamma(t) + a x^2(t),\forall t \in [0,T],
\end{align}
where $a$ is a positive constant.
Suppose that the condition $$1-4a \gamma(t)>0$$ holds for all $t \in [0,T]$ and $\xi_1(t),\xi_2(t)$ are roots of
\begin{align}
 a z^2 - z + \gamma(t) = 0,
\end{align}
$\xi_1(t)<\xi_2(t)$.
Then, if $x(0) < \xi_1(0)$, it follows that $x(t) \leq \xi_1(t),\forall t \in [0,T]$.
\end{lemma}

Using this Lemma, we assume that
\begin{align*}
& \left\{ c_0 c_1 + 2 c_1 c_2 \left(\int_{0}^{T}\|\mathbf{f}(s)\|\; ds \right) -\nu_1 \right\}^2 >\nonumber\\
&  4 c_1^2 \left\{ c_0 c_2 \int_{0}^{T}\|\mathbf{f}(s)\| ds  +  c_2^2 \left(\int_{0}^{T}\|\mathbf{f}(s)\|\; ds\right)^2 + c_0^2\right\},
\end{align*}
which simplifies to 
\begin{align}\label{conditionreg}
 & 6 c_{\mu}^2 \mathcal{E}_1(0) + 2 \nu_1 \sqrt{2} c_{\mu}\mathcal{E}_1^{1/2}(0) + 8 c_{\mu} \nu_1 \left(\int_{0}^ {T}\|\mathbf{f}\| ds\right) < \nu_1^2.
 \end{align}

Therefore $\int_{0}^{T} |\tilde{\mathcal{L}}\mathbf{h}|^2 ds$ is bounded (by a constant that does not depend on $m$), which by \eqref{ineq3} implies that there exists a positive constant $C_E$ (which also does not depend on $m$) such that $$\mathcal{E}_1(t) \leq C_E$$
for $t \in [0,T]$.
These constants depend on the initial energy $\mathcal{E}_1(0)$ and on the following norm of $\V{f}$: $\|\V{f}\|_{L^{1}(0,T;H_{0}^{1}(\Omega))}$.
\begin{remark}
We conclude that if \eqref{conditionreg} is satisfied, there exist strong solutions for the three-dimensional initial boundary-value problem. This result improves on a previous result due to Botsenyuk (\cite{Bot96a}) on the existence of strong solutions. Botsenyuk (\cite{Bot96a}) proved under his assumptions that there exists a unique strong solution. Under our assumptions the same result holds with the same proof.
\end{remark}
In the next theorem, the weak solution we are referring to are the ones proved by Botsenyuk (\cite{Bot92}).
\begin{theorem}
Let $T$ be a positive real number. Let $\mathbf{f} \in L^{1}(0,T;H_{0}^{1}(\Omega))$. We assume that the following condition holds:
\begin{align}\label{cond_nu1}
 & 6 c_{\mu}^2 \mathcal{E}_1(0) + 2 \nu_1 \sqrt{2} c_{\mu}\mathcal{E}_1^{1/2}(0) + 8 c_{\mu} \nu_1\;\|\V{f}\|_{L^{1}(0,T;H_{0}^{1}(\Omega))} < \nu_1^2,
 \end{align}
 where
$$2 \mathcal{E}_1(t) := \|\mathbf{u}^{'}(t)\|^2 + |\mathcal{L}\mathbf{u}(t)|^2 + \|\mathbf{h}(t)\|^2,$$
Then, the weak solutions $(\mathbf{u},\mathbf{h})$ have the additional regularity:\\
$\mathbf{h} \in L^{2}(0,T;D(\tilde{\mathcal{L}})) \cap L^{\infty}(0,T;H_{\sigma}^{1}(\Omega))$,\\ 
$\mathbf{u} \in L^{2}(0,T;D(\mathcal{L}))$ and\\
$\mathbf{u}^{'} \in L^{\infty}(0,T;H_{0}^{1}(\Omega))$. 
\end{theorem}

\subsection{Existence of T-periodic strong solutions}
 \parag We claim that this result of regularity holds for T-periodic solutions if the system has a linear mechanical dissipation (or nonlinear in the restricted sense discussed later), because we can show the existence of T-periodic solutions using the Poincar\'e map.
 
Using the Faedo-Galerkin  framework, we take the $L^2$ inner product of each term in the first equation of the approximate problem (with $\V{\rho}(\V{u}_m^{'})=\alpha \V{u}_m^{'}$)
$$\mathbf{u}_m^{''} + \mathcal{L} \mathbf{u}_m +\rho(\mathbf{u}_m^{'})=\curl \mathbf{h}_m \times (\mathbf{h}_m + \tilde{\mathbf{H}}) + \mathbf{f} $$
by $\mathbf{u}_m^{'}$, and the second equation, which is
$$\mathbf{h}_m^{'} + \nu_1 \tilde{\mathcal{L}} \mathbf{h}_m=\curl \left[\mathbf{u_m}^{'} \times (\mathbf{h}_m + \tilde{\mathbf{H}})\right],$$
by $\mathbf{h}_m$, we obtain the fundamental identity
\begin{align}
& \frac{d \mathcal{E}}{dt} + \left(\V{\rho}(\V{u}_m^{'}),\V{u}_m^{'} \right)_2 + \nu_1 \|\mathbf{h}_m\|^2 = \left(\mathbf{f},\mathbf{u}_m^{'}\right)_2.
\end{align}
Formally, taking the $L^2$ inner product of each term in the first equation by $\varepsilon \mathbf{u}_m$ (with $\varepsilon>0$) and adding the result to the previous identity to obtain the following inequality
\begin{align}
& \frac{d}{dt}\left\{ \mathcal{E} + \varepsilon \left(\mathbf{u}_m^{'},\mathbf{u}_m\right)_2 + \frac{\alpha \varepsilon}{2}|\mathbf{u}_m|^2\right\} + (\alpha-\varepsilon) |\mathbf{u}_m^{'}|^2\nonumber\\
& + \varepsilon \|\mathbf{u}_m\|^2+ \nu_1 \|\mathbf{h}_m\|^2 \leq |\mathbf{f}|_2 |\mathbf{u}_m^{'}|_2+ |\mathbf{f}|_2 |\mathbf{u}_m|_2+ C \|\mathbf{h}_m\|^2\;\|\mathbf{u}_m\|
\end{align}
Let $G_{\varepsilon}:=\mathcal{E} + \varepsilon \left(\mathbf{u}_m^{'},\mathbf{u}_m\right)_2 + \frac{\alpha \varepsilon}{2}|\mathbf{u}_m|^2$. Since $$\frac{1}{2} \mathcal{E} \leq G_{\varepsilon}\leq (2+\alpha)\mathcal{E}$$
if $0<\varepsilon<\min(1,C_{\Omega}^{-2})$, we infer that
\begin{align}
 \frac{d}{dt}G_{\varepsilon} + 2\varepsilon\;\mathcal{E}& \leq |\mathbf{f}|_2 |\mathbf{u}_m^{'}|_2+ C_{\Omega}^{1/2}|\mathbf{f}|_2 \|\mathbf{u}_m\|+ C \|\mathbf{h}_m\|^2\;\|\mathbf{u}_m\|\nonumber\\
 & \leq C |\mathbf{f}|_2 \sqrt{\mathcal{E}} + C \sqrt{\mathcal{E}} \|\mathbf{h}_m\|^2
\end{align}
if $0<\varepsilon<\min(1,C_{\Omega}^{-2},\alpha/2,\nu_1)$. This estimate implies that
\begin{align}
 \frac{d}{dt}\sqrt{G_{\varepsilon}}(t) + \frac{\varepsilon}{2 + \alpha}\;\sqrt{G_{\varepsilon}}(t) \leq C \left(|\mathbf{f}|_2 +\|\mathbf{h}_m\|^2 \right)
\end{align}
if $0<\varepsilon<\min(1,C_{\Omega}^{-2},\alpha/2,\nu_1)$.

Therefore,
\begin{align}
 & \sqrt{\mathcal{E}}(t)  \leq \sqrt{2} \sqrt{G_{\varepsilon}}(t)\nonumber\\
 &  \leq    C_1 \exp\left(- \frac{\varepsilon}{2 + \alpha}t \right) \left(\int_{0}^{t}e^{\frac{\varepsilon s}{2+\alpha}} |\mathbf{f}|_2 (s) ds + \int_{0}^{t}e^{\frac{\varepsilon s}{2+\alpha}} \|\mathbf{h}_m\|^2 (s) ds\right)\nonumber\\
 &+\sqrt{G_{\varepsilon}(0)} \exp\left(- \frac{\varepsilon}{2 + \alpha}t \right).\nonumber
 \end{align}
 So,
 \begin{align}
\sqrt{\mathcal{E}}(t)  & \leq C_1\;\|\mathbf{f}\|_{L^{1}(0,T;L^2(\Omega))}+\frac{C_2}{\nu_1}\mathcal{E}(0) \left(1+ \|\mathbf{f}\|_{L^{1}(0,T;L^2(\Omega))}\right)\nonumber\\
 &+\frac{C_3}{\nu_1}\|\mathbf{f}\|_{L^{1}(0,T;L^2(\Omega))}^2+ \left(\sqrt{2+\alpha}\right)\sqrt{\mathcal{E}(0)}\exp\left(-\frac{\varepsilon}{2+ \alpha}t\right),
\end{align}
where we have used the estimate (obtained from the uniform bound on $\sqrt{\mathcal{E}}$): $$\nu_1\;\int_{0}^{t}\|\mathbf{h}_{m}(s)\|^2 ds \leq\mathcal{E}(0) \left(1+ \|\mathbf{f}\|_{L^{1}(0,T;L^2(\Omega))}\right) + \|\mathbf{f}\|_{L^{1}(0,T;L^2(\Omega))}^2,\;\forall t \in [0,T].$$
We conclude that the $S(\mathbf{u}_m(0),\mathbf{u}_m^{'}(0),\mathbf{h}_m(0))=(\mathbf{u}_m(T),\mathbf{u}^{'}_m(T),\mathbf{h}_m(T))$ maps every ball of radius $R \in (R_{cr},1)$ and center  at $0$ into the same ball, where \begin{align}\label{eq_condpersol}
R_{cr}=\frac{C_1  \|\mathbf{f}\|_{L^{1}(0,T;L^2(\Omega))}+\frac{C_3}{\nu_1}\;\|\mathbf{f}\|_{L^{1}(0,T;L^2(\Omega))}^2}{1-\sqrt{2+\alpha}\exp\left(-\frac{\varepsilon}{2+\alpha}T\right)-\frac{C_2}{\nu_1}\left(1+\|\V{f}\|_{L^{1}(0,T;L^{2}(\Omega))}\right)}
\end{align}
under the condition that $R_{cr} \in (0,1)$. This map is also continuous. By the Brouwer fixed point theorem, it follows that there exists at least one solution $(\mathbf{u}_m,\mathbf{u}^{'}_m,\mathbf{h}_m)$ T-periodic. The remaining task of  passing to the limit as $m \rightarrow \infty$ was described in detail in \cite{MohOliv12} and so we omit this part of the proof. The main conclusion that we obtain is the existence of a strong T-periodic solution when the dissipation is given by $\V{\rho}(\mathbf{u}^{'})=\alpha \mathbf{u}^{'}$, under the assumption that $0<R_{cr}<1$.
 
 Notice that this result requires not only the expected restrictions on the sizes of $\nu_1$ (large) and $\mathbf{f}$ (small), but also $T$ not too small.
 \begin{theorem} (Existence of strong T-periodic solutions)\\
If $\mathbf{f} \in L^{1}(0,T;H_{0}^{1}(\Omega))$ and the conditions \eqref{cond_nu1}, $0<R_{cr}<1$ hold, then there exist strong T-periodic solutions for the magnetoelastic system \eqref{eq_1.1}, \eqref{eq_1.2}, \eqref{eq_1.3}, \eqref{eq_1.4} and \eqref{eq_1.5}.
 \end{theorem}
 
\begin{remark} The drawback of this argument of existence depending on the regularity estimate for the initial boundary-value problem and the Poincar\'e map is that we have not been able to extend this argument for the full case considered in \cite{MohOliv12}.
\end{remark}

\subsection{Stability of Strong Time-periodic Solutions}

\parag We proceed to investigate the stability of strong T-periodic solutions for the magnetoelastic system subjected now to a linear dissipation $\V{\rho}(\mathbf{v}^{'})=\alpha \mathbf{v}^{'}$. For this type of dissipation we can use the result of the previous subsection since we proved the existence of T-periodic regular solutions using the initial value problem and the Poincar\'e map.

Let $(\mathbf{v},\mathbf{b})$ denote the perturbation pair of the strong solution $(\mathbf{u},\mathbf{h})$. Then, these perturbations satisfy the following system
\begin{align}
& \mathbf{v}^{''} + \mathcal{L} \mathbf{v}+\V{\rho}(\mathbf{u}^{'}+\mathbf{v}^{'})-\V{\rho}(\mathbf{u}^{'})  =\curl \mathbf{h} \times \mathbf{b} + \curl \mathbf{b} \times \mathbf{h} + \curl \mathbf{b} \times \mathbf{b}\nonumber\\
& \mathbf{b}^{'} + \nu_1 \tilde{\mathcal{L}} \mathbf{b} = \curl[\mathbf{u}^{'} \times \mathbf{b}]+\curl[\mathbf{v}^{'} \times \mathbf{h}]+\curl[\mathbf{v}^{'} \times \mathbf{b}]\label{pgoveq}\\
& \mbox{div}\;\mathbf{b}=0\nonumber
\end{align}
under the initial conditions $\mathbf{v}(0)=\mathbf{v}_{0}$, $\mathbf{v}^{'}(0)=\mathbf{v}_{1}$ and $\mathbf{b}(0)=\mathbf{b}_{0}$.

We will look for an a priori estimate for the solutions of this system. Formally, taking the $L^{2}(\Omega)$-inner product of the first equation by $\mathbf{v}^{'}$, taking the $L^{2}(\Omega)$-inner product of the second equation by $\mathbf{b}$ and adding the results, we obtain
\begin{align}
& \frac{1}{2}\frac{d}{dt}\left\{|\mathbf{v}^{'}|_2^2 + \|\mathbf{v}\|^2 + |\mathbf{b}|_2^2\right\}+ \left(\V{\rho}(\mathbf{u}^{'}+\mathbf{v}^{'})-\V{\rho}(\mathbf{u}^{'}) ,\mathbf{v}^{'}\right)_2 +\nu_1 \|\mathbf{b}\|^2\nonumber\\
& =(\curl \mathbf{h}\times \mathbf{b},\mathbf{v}^{'})_2+ (\mathbf{u}^{'}\times \mathbf{b},\curl\mathbf{b})_2
\end{align}
 Let $2 \mathcal{E}_p(t):=|\mathbf{v}^{'}|^2 + \|\mathbf{v}\|^2 + |\mathbf{b}|^2$ denote the energy of the perturbation. Then, it follows that
 \begin{align}
 \frac{d\mathcal{E}_p}{dt}+\left(\V{\rho}(\mathbf{u}^{'}+\mathbf{v}^{'})-\V{\rho}(\mathbf{u}^{'}) ,\mathbf{v}^{'}\right)_2+\nu_1 \|\mathbf{b}\|^2 & \leq |\tilde{\mathcal{L}}\mathbf{h}|_2 \|\mathbf{b}\| |\mathbf{v}^{'}|_2\nonumber\\
 &+\|\mathbf{u}^{'}\|\; \|\mathbf{b}\|^{3/2}\; |\mathbf{b}|_{2}^{1/2}
 \end{align}
 Thus, using Young's inequality we infer that
 \begin{align}\label{eq_125}
 \frac{d\mathcal{E}_p}{dt}+\left(\V{\rho}(\mathbf{u}^{'}+\mathbf{v}^{'})-\V{\rho}(\mathbf{u}^{'}) ,\mathbf{v}^{'}\right)_2+\frac{\nu_1}{2} \|\mathbf{b}\|^2 & \leq \frac{1}{\nu_1}  |\tilde{\mathcal{L}}\mathbf{h}|_2^{2} |\mathbf{v}^{'}|_2^2+ \frac{1}{\nu_1}\; \|\mathbf{u}^{'}\|^4 |\mathbf{b}|_{2}^{2}\nonumber\\
 & \leq \frac{1}{\nu_1}  |\tilde{\mathcal{L}}\mathbf{h}|_2^{2} |\mathbf{v}^{'}|_2^2+ \frac{4}{\nu_1}\; C_E^2 |\mathbf{b}|_{2}^{2}
 \end{align}
 Assuming that $\nu_1>2 \sqrt{2}\; C_{E} C_{\Omega}^{1/2}$, it follows that
 \begin{align}
 \frac{d\mathcal{E}_p}{dt} \leq \frac{1}{\nu_1}  |\tilde{\mathcal{L}}\mathbf{h}|_2^{2} |\mathbf{v}^{'}|_2^2 \leq \frac{2}{\nu_1}  |\tilde{\mathcal{L}}\mathbf{h}|_2^{2} \;\mathcal{E}_p
 \end{align}
i.e., 
\begin{align}
\frac{d\sqrt{\mathcal{E}_p}}{dt} \leq \frac{1}{\nu_1} |\tilde{\mathcal{L}}\mathbf{h}|^2\; \sqrt{\mathcal{E}_p} 
\end{align}
which implies that
\begin{align}
\sqrt{\mathcal{E}_p} \leq \sqrt{\mathcal{E}_p(0)} \;\mbox{exp}\;\left\{\int_{0}^{t} |\tilde{\mathcal{L}}\mathbf{h}|^{2} ds\right\}\leq C_h \sqrt{\mathcal{E}_p(0)},
\end{align}
since the term $\{\int_{0}^{t}  |\tilde{\mathcal{L}}\mathbf{h}|^{2} ds\}$ is bounded by a constant as a consequence of the a priori estimates that produced the strong solution $(\mathbf{u},\mathbf{h})$.
This a priori estimate places the perturbation $(\mathbf{v},\mathbf{b})$ in the usual energy space of the weak solutions.

Now, we take the $L^{2}(\Omega)$-inner product of the first equation of \eqref{pgoveq} by $\eta \mathbf{v}$, where  $\eta\in (0,\alpha)$, to get
\begin{align*}
&\eta \frac{d}{dt}(\mathbf{v}^{'},\mathbf{v})_2+\eta\left(\V{\rho}(\mathbf{u}^{'}+\mathbf{v}^{'})-\V{\rho}(\mathbf{u}^{'}) ,\mathbf{v}\right)_2 +\eta \|\mathbf{v}\|^2 -\eta |\mathbf{v}^{'}|^2\nonumber\\ 
&\leq 2\eta\|\mathbf{v}\|\; \|\mathbf{b}\|\; \|\mathbf{h}\|+\eta\|\mathbf{v}\|\;\|\mathbf{b}\| \|\mathbf{b}\|_{H^{1/2}}
\end{align*}
Adding each member of this equation to the corresponding members of equation \eqref{eq_125}, we obtain the estimate
 \begin{align}
& \frac{d}{dt}\{\mathcal{E}_p  +\eta (\mathbf{v}^{'},\mathbf{v})_{2}\}+ \eta\left(\V{\rho}(\mathbf{u}^{'}+\mathbf{v}^{'})-\V{\rho}(\mathbf{u}^{'}) ,\mathbf{v}\right)_2 +\left(\V{\rho}(\mathbf{u}^{'}+\mathbf{v}^{'})-\V{\rho}(\mathbf{u}^{'}) ,\mathbf{v}^{'}\right)_2\nonumber\\
& \quad \quad \quad \quad \quad \quad \quad \quad -\eta |\mathbf{v}^{'}|^2+\eta \|\mathbf{v}\|^2+\frac{\nu_1}{2} \|\mathbf{b}\|^2 \nonumber\\
 & \leq  \frac{1}{\nu_1}\; |\tilde{\mathcal{L}}\mathbf{h}|_2^{2} |\mathbf{v}^{'}|_2^2+ c\; \|\mathbf{u}^{'}\|^4 |\mathbf{b}|_{2}^{2}+ 2\eta\|\mathbf{v}\|\; \|\mathbf{b}\|\; \|\mathbf{h}\|+\eta\|\mathbf{v}\|\;\|\mathbf{b}\|^{3/2}\;|\mathbf{b}|_2^{1/2}
  \end{align}
We infer that
 \begin{align}
& \frac{d}{dt}\{\mathcal{E}_p+\eta (\mathbf{v}^{'},\mathbf{v})_{2}+\frac{\eta \alpha}{2}|\mathbf{v}|_2^2\} +\left(\alpha -\eta\right) |\mathbf{v}^{'}|^2+\eta \|\mathbf{v}\|^2+\frac{\nu_1}{2} \|\mathbf{b}\|^2 \nonumber\\
 & \leq  \frac{1}{\nu_1}\; |\tilde{\mathcal{L}}\mathbf{h}|_2^{2} |\mathbf{v}^{'}|_2^2+ c\; \|\mathbf{u}^{'}\|^4 |\mathbf{b}|_{2}^{2}+ 2\eta\|\mathbf{v}\|\; \|\mathbf{b}\|\; \|\mathbf{h}\|+\eta\|\mathbf{v}\|\;\|\mathbf{b}\|^{3/2}\;|\mathbf{b}|_2^{1/2}
  \end{align}
  Let $$G_p=\mathcal{E}_p+\eta (\mathbf{v}^{'},\mathbf{v})_{2}+\frac{\eta \alpha}{2}|\mathbf{v}|_2^2.$$
  The next lemma establishes the relation between the function $G_p$ and the energy.
\begin{lemma} If $0<\eta<\min(1,(C_{\Omega})^{-2})$, then
\begin{align}\label{eq_GE}
\frac{1}{2}\mathcal{E}_p(t)\leq G_p(t) \leq (2+\alpha) \mathcal{E}_p(t).
\end{align}
\end{lemma}
\begin{proof}
\begin{align}\label{eq_412}
\frac{1}{2}\mathcal{E}_p &\leq G_p \leq \mathcal{E}_p +
\frac{1}{2}|\V{v}^{\prime}|_{2}^{2}+ |\V{v}|_{2}^{2}\frac{\eta^2
  C_{\Omega}^2}{2 C_{\Omega}^2} +\frac{\eta \alpha}{2}|\mathbf{v}|_2^2\nonumber\\
& \leq |\V{v}^{\prime}|_2^2 + \frac{1}{2}|\V{b}|_2^2 + \left(
\frac{1}{2}+\frac{\eta}{2 (C_{\Omega})^{-2}}+\frac{\alpha \eta}{2 C_{\Omega}^{-2}}\right) a_{I}(\V{v},\V{v})
\leq (2+\alpha) \mathcal{E}_p. \nonumber
\end{align}
\end{proof}

  Then,
\begin{align}  
& \frac{dG_p}{dt} +\left(\alpha -\eta\right) |\mathbf{v}^{'}|^2+\eta \|\mathbf{v}\|^2+\frac{\nu_1}{2} \|\mathbf{b}\|^2 \nonumber\\
 & \leq  \frac{1}{\nu_1}\; |\tilde{\mathcal{L}}\mathbf{h}|_2^{2} |\mathbf{v}^{'}|_2^2+ c\; \|\mathbf{u}^{'}\|^4 |\mathbf{b}|_{2}^{2}+ 2\eta\|\mathbf{v}\|\; \|\mathbf{b}\|\; \|\mathbf{h}\|+\eta\|\mathbf{v}\|\;\|\mathbf{b}\|^{3/2}\;|\mathbf{b}|_2^{1/2}
  \end{align}
So, for $\eta \in (0,\alpha/2)$ sufficiently small,
 \begin{align}
& \frac{d G_p}{dt}+\frac{\alpha}{2}|\mathbf{v}^{'}|_2^2+\frac{\eta}{2} \|\mathbf{v}\|^2 +\left\{\frac{C_{\Omega}^{-1/2}\nu_1}{4}-c\;C_E^2\right\} |\mathbf{b}|_2^2\nonumber\\
&+\left\{\frac{\nu_1}{4}-2\eta\;C_E\right\} \|\mathbf{b}\|^2 \leq \eta\|\mathbf{v}\|\;\|\mathbf{b}\|^{3/2}\;|\mathbf{b}|_2^{1/2}+ \frac{1}{\nu_1}\;|\tilde{\mathcal{L}}\mathbf{h}|_2^{2}\; |\mathbf{v}^{'}|_2^2.
  \end{align}
 This inequality and the previous uniform bound on the energy of perturbations imply the following estimate
 \begin{align}
& \frac{d G_p}{dt}+\frac{\alpha}{2} |\mathbf{v}^{'}|_2^2+\frac{\eta}{2} \|\mathbf{v}\|^2 +\left\{\frac{C_{\Omega}^{-1/2}\nu_1}{4}-c\;C_E^2\right\} |\mathbf{b}|_2^2\nonumber\\
&+\left\{\frac{\nu_1}{4}-2\eta\;C_E \right\} \|\mathbf{b}\|^2  \leq  C_h \sqrt{2 \mathcal{E}_p(0)}\; \eta\; \mathcal{E}_p + \frac{2}{\nu_1}\;|\tilde{\mathcal{L}}\mathbf{h}|_2^2\; \mathcal{E}_p
  \end{align}
Thus, 
 \begin{align}
& \frac{d G_p}{dt} + \left(C_0 -C_h \sqrt{2 \mathcal{E}_p(0)}\; \right)\eta\; \mathcal{E}_p \leq \frac{2}{\nu_1}\;|\tilde{\mathcal{L}}\mathbf{h}|_2^2\; \mathcal{E}_p  \end{align}
 Now, using the previous lemma,
 \begin{align}
& \frac{d G_p}{dt} + \frac{C_0 -C_h \sqrt{2 \mathcal{E}_p(0)}}{(2+\alpha)}\; \eta\;G_p \leq  \frac{4}{\nu_1}\;|\tilde{\mathcal{L}}\mathbf{h}|_2^2 G_p(t)
  \end{align}
  Let $$C_1:=\frac{C_0 -C_h \sqrt{2 \mathcal{E}_p(0)}}{(2+\alpha)}$$ and $$C_h:=\text{sup}\;\left\{\int_{0}^{t}|\tilde{\mathcal{L}}\mathbf{h}|_2^2 dt, t>0\right\},$$
which exists by the previous estimates.
Then, we infer that
$$G_p(t) \leq G_p(0) \exp\left(-C_1\; \eta\; t+\frac{4C_h}{\nu_1}\right)$$
 This inequality implies, for $$0 < \eta < \text{min}\;\left\{1,\frac{\alpha}{2 },\frac{1}{C_{\Omega}^{2}}\right\}$$ and 
 \begin{align}\label{cond2_nu1}
 \nu_1>\max\left\{\sqrt{2} \;C_E,2c\;C_E^2\right\}\;2C_{\Omega}^{1/2},\quad \mathcal{E}_p(0) < \frac{C_{0}^2}{2\;C_h^2}
 \end{align}
 the following main result
 \begin{theorem} (Conditional Asymptotic Stability)\\
 If the conditions \eqref{cond_nu1}, $0<R_{cr}<1$ and \eqref{cond2_nu1} hold, then\\
 (i) there exist weak solutions for the initial-value problem \eqref{pgoveq} for the perturbation field $(\mathbf{v},\mathbf{b})$,\\
 (ii) the energy of the perturbations has an exponential decay to zero as $t \rightarrow \infty$:
$$\mathcal{E}_p(t) \leq 2\;(2+\alpha)\;\mathcal{E}_p(0) \exp\left(-C_1\;\eta\; t+\frac{4C_h}{\nu_1}\right).$$
 The same result holds for the case with nonlinear dissipation $\V{\rho}(\V{u}^{'})$ (with $K_c$ replacing $\alpha$) under the assumption that there exists
a positive constant $K_c$  such that $$\left(\V{\rho}(u+w)-\V{\rho}(u),z\right)_E \geq K_c \left(w,z\right)_E, \forall u,w,z \in \mathbb{R}^3$$ 
 \end{theorem}
 \section{PART II: Asymptotic Stability for a Magnetoelastic System in 2D}
\subsection{Introduction}
\hspace{0.65cm} We consider now the asymptotic stability of solutions of an initial boundary value problem that models the interaction
between a conducting nonferromagnetic homogeneous isotropic elastic body that occupies at an initial time a bounded, simply-connected domain
$\Omega \subset \mathbb{R}^2$ with smooth boundary. We suppose that a time-independent magnetic field $(0,0,B_0)$ with 
$B_0:\Omega \rightarrow \mathbb{R}$ affects this region. Let $\mathbf{u} = (u_1,u_2,0)$ denote the displacement field and $(0,0,h(x_1,x_2))$
denote a disturbance of the magnetic field. Let $T$ be a fixed positive number. We also denote by $\mathbf{u}^{'}$ the partial derivative of
$\mathbf{u}(x,t)$ with respect to $t$ (time) and $\mathbf{u}^{''}$ means $\left(\mathbf{u}^{'}\right)^{'}$. Analogous notation is used for the scalar fields.

The corresponding system of partial differental equations is the following:
\begin{align}
&\rho_M \frac{\partial^2 \mathbf{u}}{\partial t^2} - \mu\Delta \mathbf{u} - (\lambda+\mu)\nabla \mbox{div}\;\mathbf{u}
+ \mu_0\; (B_0+h)\;\nabla h = \mathbf{f}_2 \label{eq_II_1.1}\\
& \frac{\partial h}{\partial t} - \nu_1\;\Delta h + \mbox{div}\;\left((B_0+h)\;\frac{\partial \mathbf{u}}{\partial t}\right) = f_1 \label{eq_II_1.2}
\end{align}
in $Q_T=\Omega \times (0, T)$. $\rho_M$ is the mass density per unity volume of the medium. $\lambda$ and $\mu$ are Lam\'e's constants $(\mu > 0$,
$\lambda > 0)$. $\nu_1 = 1/(\sigma \mu_0)$ where $\sigma>0$ represents the conductivity of the material and 
$\mu_0$ is a positive number representing the magnetic permeability. $\mathbf{f}=(f_1,\mathbf{f}_2)$ is a known external force. The initial
and boundary conditions associated with the system are
\begin{equation}\label{eq_II_1.3}
\mathbf{u}(x,0) = \mathbf{u}_0(x), \quad \frac{\partial \mathbf{u}}{\partial t}(x,0) = \mathbf{u}_1(x),\quad h(x,0) = h_0(x)\;\mbox{in}\; \Omega
\end{equation}
and 
\begin{equation}\label{eq_II_1.4}
\mathbf{u}=0, \quad \nabla b \cdot \mathbf{n} = 0 \;\;\mbox{on}\;\;\Sigma_T:=\partial\Omega \times (0,T),
\end{equation}
where $\mathbf{n}=(n_1,n_2) = \mathbf{n}(x)$ denotes the outer unit normal at $x=(x_1,x_2) \in \partial\Omega$.

The total energy associated with the above system is given by
\begin{equation}
\mathcal E(t) = \frac{1}{2}\;\left\{\rho_M \left|\mathbf{u}^{'}\right|^2_{L^2} + \mu|\nabla \mathbf{u}|^2_{L^2}
+ (\lambda+\mu)|\mbox{div}\; \mathbf{u}|^2_{L^2} + |h|^2_{L^2}\right\} \label{eq_II_1.5}
\end{equation}
where $|\;|_{L^2}$ and $(\;,\;)_{L^2}$ denote the $L^2(\Omega)$-norm and related inner product, respectively.

If $f=0$, then, formally, by taking the $L^2$-inner product of each term in equation (\ref{eq_II_1.1}) by $\mathbf{u}^{'}$ and each term in 
(\ref{eq_II_1.2}) by $b$ and adding the resulting relations, we obtain the \emph{energy identity}
\begin{equation}
\frac{d}{dt}\,\mathcal{E}(t) = -\nu_1 |\nabla h|^2 \label{eq_II_1.6}
\end{equation}
which means that the model (\ref{eq_II_1.1})-(\ref{eq_II_1.5}) has a natural dissipative term. This identity also provides the first a priori estimate (\cite{Bot96b}) for each $t \in [0,T)$:
\begin{equation}
|h(t)|^2_{L^2}+\int_{0}^{t}|\nabla h(t)|^2_{L^2}\;ds + |\mathbf{u}^{'}(t)|^2_{L^2} + \|\mathbf{u}(t)\|^2 \leq c_1 \label{eq_II_1.7}
\end{equation}

Botsenyuk (\cite{Bot96b}) proved the existence and uniqueness of global solutions for the two-dimensional problem without
assuming any restrictions on the initial data. The following estimate was obtained:
\begin{equation}\label{eq_II_3.1}
\int_{0}^{t} \|h\|_{1+\gamma}^2 ds + \|h\|_{\gamma}^2 +  \|\mathbf{u}^{'}\|_{\gamma}^2 +  \|\mathbf{u}\|_{1+\gamma}^2  \leq c_2,
\end{equation}
where $\|\;\|_{s}$ denotes the norm of $H^s(\Omega)$, $s \in \mathbb{R}$, and $0<\gamma<1/2$.

We investigate the large time behavior of the solutions of this system without any additional damping
mechanism (and no external force).

\subsection{Weak formulation}
Let $T>0$, $B_0 \in L^2(0,T;H^1(\Omega)\cap L^{\infty}(\Omega))$, $h_0 \in H^{\gamma}(\Omega)$, $\mathbf{u}_0 \in H_0^{1+\gamma}(\Omega)$, 
$\mathbf{u}_1 \in H^{\gamma}(\Omega)$, $f_1 \in L^{2}(0,T;L^{2}(\Omega))$ and $\mathbf{f}_2 \in L^{1}(0,T;H^{\gamma}(\Omega))$ $(0 < \gamma < 1)$. We look for
functions $(\mathbf{u},\mathbf{u}^{'},h)$ in $$L^{\infty}(0,T;H_0^1(\Omega))\times L^{\infty}(0,T;L^2(\Omega)) \times [L^2(0,T;H^1(\Omega)\cap L^{\infty}(0,T;L^2(\Omega))]$$
satisfying the following system:
\begin{align}
& \rho_M \frac{d}{dt}(\mathbf{u}^{'},\psi)+a_2(\mathbf{u},\psi)+\mu_0\;\int_{\Omega} (B_0+h)\;\nabla h\cdot \psi\;dx=(\mathbf{f}_2,\psi),\label{eq_2.1B}\\
& \forall \psi \in H_0^{1}(\Omega)\nonumber\\
& \frac{d}{dt}(h,\varphi)+a_1(h,\varphi)- (h,\varphi) -\int_{\Omega} (B_0+h)\; \mathbf{u}^{'}\cdot \nabla\varphi\;dx=(f_1,\varphi),\label{eq_2.2B}\\
& \forall \varphi \in H^{1}(\Omega)\nonumber
\end{align}
where $$a_1(\varphi,\eta):=\nu_1\;(\nabla \varphi,\nabla \eta)+ (\varphi, \eta),\;\forall \varphi,\eta \in H^1(\Omega)\quad \mbox{and}$$
$$a_2(\mathbf{v},\mathbf{w}):=\mu\;(\nabla \mathbf{v}, \nabla \mathbf{w}) + (\lambda+\mu) (\;\mbox{div}\;\mathbf{v},\;\mbox{div}\;\mathbf{w}),\;\forall \mathbf{v},\mathbf{w} \in H_0^1(\Omega).$$

We will make use of the following existence-uniqueness result of Botsenyuk:
\begin{theorem}(Botsenyuk, \cite{Bot96b})\\
Under the previous assumption on the initial data and forcing function, then there exist unique functions $(\mathbf{u},h)$ satisfying \eqref{eq_1.3}, \eqref{eq_1.4}, \eqref{eq_2.1B} and \eqref{eq_2.2B}. Moreover,
$h \in L^2(0,T;H^{1+\gamma}(\Omega)),$ $h \in C([0,T];H^{\gamma}(\Omega))$, $\mathbf{u} \in C([0,T];H_{0}^{\gamma+1}(\Omega))$ and
$\mathbf{u}^{'} \in C([0,T];H^{\gamma}(\Omega))$.
\end{theorem}

\subsection{Asymptotic behavior}
We shall assume that $B_0$ is a constant and $\mathbf{f} = 0$. We denote by $L^2_0(\Omega)$ the set of functions $\varphi$ in $L^{2}(\Omega)$ that satisfy $\int_{\Omega} \varphi \;dx=0$.
We must have $h_0 \in H^{\gamma}(\Omega) \cap L^2_0(\Omega)$ and the solution $h$ belongs to the spaces $L^2(0,T;H^{1+\gamma}(\Omega)\cap L^2_0(\Omega))$ and $C([0,T];H^{\gamma}(\Omega)\cap L^2_0(\Omega))$,
since (considering the three-dimensional problem) the identity $\partial_j(x_k h_j)=(\partial_j x_k)h_j$, integration in $\Omega$, the boundary condition $\mathbf{h}\cdot n=0$ and the divergence theorem imply that $h$ must satisfy 
\begin{equation}\label{hmeanzero}
\int_{\Omega} h\;dx=0.
\end{equation}

We recall some definitions and results concerning dynamical systems.
\begin{definition}
A dynamical system on a Banach space $Z$ is a family\\ $\{S(t)\}_{t \geq 0}$ of mappings on $Z$ such that\\
(i) $S(t) \in C(Z,Z),\forall t \geq 0$,\\
(ii) $S(0)=I$,\\
(iii) $S(t+s)=S(t) S(s),\forall s,t \in \mathbb{R}_{+}$,\\
(iv) the function $t \mapsto S(t)z$ belongs to $C([0,\infty),Z)$,$\forall z \in Z$.\\
An orbit (positive orbit) $\gamma^{+}=\gamma^{+}(\varphi)$ through $\varphi$ in $Z$ is defined to be $\gamma^{+}(\varphi)=\cup_{t \geq 0}S(t)\varphi$.
\end{definition}
\begin{definition}
Let $\{S(t)\}_{t \geq 0}$ is a dynamical system on a Banach space $Z$ and $V$ is a continuous scalar function defined on $Z$.
The function $\dot{V}(\varphi)$ is defined by $$\dot{V}(\varphi)=\lim \sup_{t \rightarrow 0^{+}}\frac{1}{t}\left[V(S(t)\varphi)-V(\varphi)\right].$$
$V:Z \rightarrow \mathbb{R}$ is a Lyapunov function on a set $G \subset Z$ if $V$ is continuous on the closure of $G$ and $\dot{V}(\varphi)\leq 0$ for
$\varphi \in G$.
\end{definition}

\begin{theorem}(\cite{Hal69})\\
Let $\{S(t)\}_{t \geq 0}$ is a dynamical system on a Banach space $Z$. If $V$ is a Lyapunov function on a set $G \subset Z$ and an
orbit $\gamma^{+}(z)$ belongs to $G$ and is in a compact set of $Z$, then $S(t)(z) \rightarrow M$ as $t \rightarrow \infty$, where $M$ is the
largest invariant set in $S:=\{\varphi \in \overline{G}:\dot{V}(\varphi)=0\}$ of $\{S(t)\}$.
\end{theorem}

Now, let $X=H_{0}^{1+\gamma}(\Omega) \times H^{\gamma}(\Omega) \times \left[H^{\gamma}(\Omega)\cap L^2_0(\Omega)\right]$ with $0<\gamma < 1/2$. Let
$Y=H_{0}^{1+\gamma^{'}}(\Omega) \times H^{\gamma^{'}}(\Omega) \times \left[H^{\gamma^{'}}(\Omega)\cap L^2_0(\Omega)\right]$ with $0<\gamma^{'}<\gamma < 1/2$. The total energy of the magnetoelastic system is defined by
$$\mathcal{E}(t)=\frac{1}{2}\int_{\Omega} \left[\rho_M |\mathbf{u}_t(t)\right|^2 + \mu |\nabla \mathbf{u}(t)|^2+(\lambda + \mu)\;|\mbox{div}\;\mathbf{u}(t)|^2+h^2(t)]\;dx.$$

We know that $$\frac{d\mathcal{E}}{dt}=-\frac{1}{\sigma \mu_e}\int_{\Omega}|\nabla h|^2\;dx\leq 0,$$
which implies
\begin{equation}\label{eq_II_3p1}
\mathcal{E}(t) +\frac{1}{\sigma \mu_e}\int_{0}^{t}\int_{\Omega}|\nabla h|^2\;dx\;ds=\mathcal{E}(0), \quad \forall t \geq 0.
\end{equation}
Let
\begin{equation}
\mathcal{E}_{\gamma}(t):=\|\mathbf{u}_t(t)\|_{\gamma}^2 + \|\mathbf{u}(t)\|_{1+\gamma}^2+\|h(t)\|_{\gamma}^2.    
\end{equation} 
We also know (\cite{Bot96b}) that there exists a positive constant $C_2$, which depends on $\mathcal{E}(0)$ and also on $\|\mathbf{u}_1\|_{\gamma}^2$,$\|\mathbf{u}_0\|_{1+\gamma}^2$,
$\|h_0\|_{\gamma}^2 $, such that 
\begin{equation}\label{eq_egamma}
\mathcal{E}_{\gamma}(t) + \int_{0}^{t}\|h\|_{\gamma+1}^2\;ds\leq C_2,\quad \forall t \geq 0.
\end{equation}
It follows from the well-posedness result of Botsenyuk (\cite{Bot96b}) that $\{S(t)\}_{t \geq 0}$ given by $S(t):X \rightarrow Y$, $S(t)(\mathbf{u}_0,\mathbf{u}_1,h_0)=(\mathbf{u},\mathbf{u}_t,h)$ ($t \geq 0$),
is a dynamical system.

The continuous function $V: X \rightarrow \mathbb{R}$ defined by $V(z)=\mathcal{E}(S(t)(z)),\;\forall z \in X$ is a Lyapunov function for $\{S(t)\}_{t \geq 0}$.

It follows from \eqref{eq_II_3p1}, \eqref{eq_egamma} that for any initial data in $X$, the trajectory remains in a bounded set of $X$, for any $t > 0$. By compactness of the Sobolev embeddings,
the trajectory remains in a compact set of $Y$. The hypothesis of Hale's theorem (\cite{Hal69}) are satisfied choosing $Z=Y$, $G=X$. Here, $M$ is the largest (positively) invariant set in
$$S=\{(\mathbf{v},\mathbf{w},b) \in X:\int_{\Omega}|\nabla b|^2\;dx=0\}.$$ 
So,
$$M=\{(\mathbf{v},\mathbf{w},b) \in X: (\mathbf{u}(t),\mathbf{u}^{'}(t),h(t))=S(t)(\mathbf{v},\mathbf{w},b) \in X,$$
$$\int_{\Omega}|\nabla h(t)|^2\;dx=0,\;\forall t \geq 0\},$$
i.e., $$M=\{(\mathbf{v},\mathbf{w},b) \in X: (\mathbf{u}(t),\mathbf{u}^{'}(t),h(t))=S(t) (\mathbf{v},\mathbf{w},b) \in X,$$
$$h(x,s)=h(s)\; \mbox{a.e.}\;\mbox{in}\; \Omega, \forall t \geq 0\}.$$
Choosing $\varphi \equiv 1$ in \eqref{eq_2.2B}, we conclude that $h=\mbox{constant}$. It follows from \eqref{hmeanzero}, that $h=0$.\\
Since we are assuming that $B_0$ is constant, it follows that $M$ is given by
 $$M=\{(\mathbf{v},\mathbf{w},b) \in X: (\mathbf{u}(t),\mathbf{u}^{'}(t),h(t))=S(t) (\mathbf{v},\mathbf{w},b) \in X,$$
 $$h(x,s)=0,\;\mbox{div}\;\mathbf{u}^{'}=0,\;\rho_M \mathbf{u}^{''}+\mathcal{L}\mathbf{u}=0\;\mbox{in}\; \Omega \times[0,\infty)\}.$$
If, in addition, $\Omega$ has the following property $\mathcal{P}$ (\cite{LeZu99})
\begin{align}\label{eq_prop}
&\mbox{If}\;\xi \in H_{0}^{1}(\Omega) \mbox{is such that}\nonumber\\
&(i) -\Delta \xi = \gamma^2 \xi,\; (ii)\mbox{div}\;\xi=0\quad \mbox{in}\; \Omega,\nonumber\\
&(iii)\xi=0\quad\mbox{on}\;\partial \Omega\\
&\mbox{for some}\;\gamma \in \mathbb{R},\nonumber\\
&\mbox{imply that}\;\xi \equiv 0,\nonumber
\end{align}
then, $M=\left\{(0,0,0)\right\}.$

We have proved the following result on the asymptotic behavior of the system under the action of only the natural dissipation that comes from the equation for the magnetic field:
\begin{theorem}
Let $\Omega$ be a bounded, simply connected domain in $\mathbb{R}^{2}$ with boundary of class $C^2$. Let $0 < \gamma_1 <\gamma_2<1/2$. Let $\{S(t)\}_{t \geq 0}$
be the dynamical system defined on $X_{\gamma_1}=H_{0}^{1+\gamma_1}(\Omega) \times H^{\gamma_1}(\Omega) \times \left[H^{\gamma_1}(\Omega)\cap L^2_0(\Omega)\right]$ by
$S(t)(\mathbf{u}_0,\mathbf{u}_1,h_0)=(\mathbf{u}(t),\mathbf{u}_t(t),h(t))$ ($t \geq 0$). We assume that $B_0$ is constant. Then,\\
(i) If $\Omega$ has the property $\mathcal{P}$, then
\begin{align}
& \mathbf{u} \rightarrow 0\; (\mbox{in}\;H_{0}^{1+\gamma_1}(\Omega)),\;h \rightarrow 0\;(\mbox{in}\;H^{\gamma_1}(\Omega)\cap L^2_{0}(\Omega)),\;\mbox{as}\; t \rightarrow \infty.\nonumber
\end{align}
(ii) If $\Omega$ does not satisfy property $\mathcal{P}$ (e.g., $\Omega$ is a disk), then for any initial data
in $X_{\gamma_2}$, $(\mathbf{u},\mathbf{u}_t,h)(t)$ aproaches the set $M$ in the norm of the space
$X_{\gamma_1}$ as $t \rightarrow \infty$, i.e., 
\begin{align}
& \mathbf{u} \rightarrow \mathbf{v}\; (\mbox{in}\;H_{0}^{1+\gamma_1}(\Omega)),\;h \rightarrow 0\;(\mbox{in}\;H^{\gamma_1}(\Omega)\cap L^2_{0}(\Omega)),\;\mbox{as}\; t \rightarrow \infty,\nonumber
\end{align}
where $\mathbf{v}$ satisfies $\mbox{div}\;\mathbf{v}^{'}=0,\;\mathbf{v}^{''}+\mathcal{L} \mathbf{v}=0$, where
$\mathcal{L}\mathbf{u}:=- \mu\Delta \mathbf{u} - (\lambda+\mu)\nabla \mbox{div}\;\mathbf{u}$.
\end{theorem}
\begin{remark}
Related to the property $\mathcal{P}$ are the following facts:\\
(i) Weinacht (\cite{Wei99}) proved that if $\Omega$ is a simply connected two-dimensional bounded region with boundary of class $C^2$ and the eigenvalue problem
\begin{align}
\mathcal{L} \mathbf{z}_l=\lambda_l\; \mathbf{z}_l\;\mbox{in}\;\Omega,\\
\mbox{div}\;\mathbf{z}_l=0\;\mbox{in}\;\Omega,\\
\mathbf{z}_l|_{\partial \Omega}=0
\end{align}
has classical non-trivial solutions for each $l=1,2,\cdots$, then $\Omega$ is a disk.\\
(ii) In Dafermos (\cite{Daf68}), it is observed that for a disk $D=\{x \in \mathbb{R}^2:|x|_E \leq 1\}$, there are non-zero eigenfunctions $v^m=(v^m_{r},v^m_{\theta}=(0,J_1(\zeta_m r))$,
$m \in \mathbb{N}$, where $J_1$ is the Bessel functions of the first kind of order 1 and $\zeta_m$ is the m-th non-negative root of $J_1$.
\end{remark}
\bibliographystyle{spmpsci}      
\bibliography{stabilityME-jco}
\end{document}